\documentclass[12pt]{amsart}
\usepackage{amssymb,amsmath,latexsym,enumerate,graphicx,bbm,mathptmx,microtype}
\usepackage{amsfonts}
\usepackage{epstopdf}
\usepackage{verbatim}
\usepackage{hyperref}
\usepackage[latin1]{inputenc}
\usepackage[english]{babel}
\usepackage[mathscr]{eucal}
\usepackage{stmaryrd}
\usepackage[arrow,matrix]{xy}
\usepackage{enumerate}
\makeatletter
\@namedef{subjclassname@2010}{
  \textup{2010} Mathematics Subject Classification}
\makeatother

\newtheorem{thm}{Theorem}[section]

\newtheorem{lem}[thm]{Lemma}

\def\q{\quad}

\def\t{\hbox}
\def\mod#1{\ (\hbox{\rm mod}\ #1)}

\def\qtq#1{\q\t{#1}\q}
\def\mod#1{\ (\text{\rm mod}\ #1)}
\def\qtq#1{\q\t{#1}\q}
\def\f{\frac}
\def\e{\equiv}
\def\b{\binom}

\def\sls#1#2{(\f{#1}{#2})}

\theoremstyle{definition}

\theoremstyle{remark}

\makeatletter
\AtBeginDocument{\@ifpackageloaded{natbib}{\ifNAT@numbers\if@filesw\immediate\write\@auxout{\string\global\string\NAT@numberstrue}\fi\fi}{}}
\makeatother
\begin{document}

\title {On the properties
of a sequence concerning binomial coefficients}

\author{Daeyeoul Kim${}^{(1)}$}
\address{${}^{(1)}$ National Institute for Mathematical Sciences
 \\ Yuseong-daero 1689-gil \\ Yuseong-gu \\
Daejeon 305-811 \\ South Korea} \email{daeyeoul@nims.re.kr}

\author{Ayyadurai Sankaranarayanan${}^{(2)}$}
\address{${}^{(2)}$
  School of
Mathematics, Tata
Institute of Fundamental Research\\Homi Bhabha Road \\Mumbai 400005\\
India}  \email{sank@math.tifr.res.in}

\author{Zhi-Hong Sun${}^{(3)}$}

\address{${}^{(3)}$School of Mathematical Sciences\\
Huaiyin Normal University\\
Huaian, Jiangsu 223001\\
The People's Republic of China} \email{zhihongsun@yahoo.com}
\thanks{\footnotesize{\it Corresponding Author:}
Zhi-Hong Sun (zhihongsun@yahoo.com)}

\begin{abstract}  For $n\ge 3$ let $f(n)$  be the least positive integer $k$
such that $\binom nk>\frac{2^n}{n+1}$. In this paper we investigate
the properties of $f(n)$. \par\q\newline MSC(2010): 11B65, 11B83,
05A20
\newline Keywords: sequence, binomial coefficient, inequality
\end{abstract}

\maketitle

\section{\bf Introduction }
\par
Throughout this paper, $n$ denotes a positive integer,
 $[x]$  denotes the greatest integer
not exceeding $x$, and  $\lceil x \rceil$ denotes  the least integer
greater than or equal to $x$. For $n\geq 3$ and $k\in\{ 1,2,\cdots,
n\}$, we have $$ \binom{n}{k} \leq \binom{n}{ [\frac{n}{2}] }
\qtq{and hence}(n-1)\binom{n}{[\frac{n}{2}]}\ge
\sum_{k=1}^{n-1}\binom{ n}{k} =2^n-2.$$ Using the fact that $2^{n-1}
> n$ for $n\ge 3$, we obtain
\begin{align*}\label{1.1}
 \binom{n}{ [ \frac{n}{2}] } \geq \frac{2^n -2}{n-1}
> \frac{2^n}{n}> \frac{2^n}{n+1}.\tag {1.1}
\end{align*}

In view of the above, for  $n\geq 3$ we define $f(n)$ to be the
least positive integer $k$ such that $\b nk>\f{2^n}{n+1}$, and
define $L(n)$ to be the least positive integer $k$ such that
$\binom{ n}{k}
>\frac{2^n}{n}$. Then  clearly
\begin{align}\label{1.2}
f(n)\le L(n) \leq [\frac{n}{2}].\tag {1.2}
\end{align} The first twenty-one values of $f(n)$ and
$L(n)$  are given below :

\begin{center}
\begin{tabular}{|c|c|c|c|c|c|c|c|c|c|c|c|c|c|c|c|c|c|c|c|c|c|}
 \hline

 $n$ & $3$ & $4$ & $5$ & $6$ & $7$ & $8$ &$9$ & $10$ & $11$ & $12$ & $13$ & $14$
 & $15$ & $16$ & $17$ & $18$ & $19$ & $20$ & $21$ & $22$ & $23$ \\
\hline $f(n)$ & $1$ & $1$ & $2$ & $2$ & $2$ & $3$ &$3$ & $3$ & $4$ &
$4$ &
 $4$ & $5$ & $5$ & $5$ & $6$ & $6$ & $6$ & $7$ & $7$ & $8$ & $8$ \\
\hline
$L(n)$ & $1$ & $2$ & $2$ & $2$ & $2$ & $3$ &$3$ & $3$ & $4$ & $4$ & $4$ & $5$ & $5$ & $5$ & $6$ & $6$ & $7$ & $7$ & $7$ & $8$ & $8$ \\
\hline
\end{tabular}
\end{center}
\begin{center}   Table 1.  Values of $f(n)$ and $L(n)$
\end{center}
The goal of this paper is to study certain properties of the
function $f(n)$. We establish the following results:

\begin{thm}\label{thm4} For $n\geq 3$ we have $f(n+1)\in\{f(n),f(n)+1\}$ and
$L(n+1)\in\{L(n),L(n)+1\}$.
\end{thm}

\begin{thm}\label{thm1-1} For $n\ge 4$ we have
$f(n)\in\{f(n-1),f(n+1)\}$ and $L(n)\in\{L(n-1),L(n+1)\}$.
\end{thm}
\begin{thm}\label{thm1-1}
For $n\ge 3$ we have $f(n)=L(n)$ or $L(n)-1$.
\end{thm}
\noindent{\bf Corollary 1.1}  Let $n\ge 3$. Then
$$f(n)\not=L(n)\iff f(n)=L(n)-1\iff \b n{L(n)-1}>\f{2^n}{n+1}.$$

\begin{thm}\label{Theorem 1.4}
Suppose $n\ge 3$ and $\b n{L(n)-1}>\f{2^n}{n+1}$. Then
$$f(n)=f(n-1)=f(n-2)=L(n)-1\qtq{and}L(n)=L(n+1)=L(n+2)=f(n)+1.$$
\end{thm}
\begin{thm}\label{Theorem 1.5} For $n\ge 3$ we have $f(n+3)>f(n)$.
\end{thm}
\noindent{\bf Remark 1.1} When $f(n)$ is replaced by $L(n)$, Theorem
1.5 is also true. Using Stirling's formula for the Gamma function
$\Gamma(x)$ (see [1-3]), one may deduce  that
$$f(n)=\f n2-\f 12\sqrt{n\hbox{log}\;\f{2n}{\pi}}+\hbox{O}(1).$$

 \section{\bf Two lemmas}
\begin{lem}\label{Lemma 2.1} Let $k\in\{0,1,\ldots,n-1\}$. Then
$$\b nk\le \f{n^n}{k^k(n-k)^{n-k}}.$$
\end{lem}
Proof. Clearly the inequality holds for $k=0$. Now we assume $1\le
k<n$. Let
$$f(x)=(1+x)^{\f 1x}\qtq{and} F(x)=\log f(x)=\f 1x\log (1+x).$$
 Then
 $$F'(x)=-\f 1{x^2}\log (1+x)+\f 1x\cdot\f 1{1+x}=\f 1{x^2}
 \big(\f x{1+x}-\log (1+x)\big).$$
Set $g(x)=\f x{1+x}-\log (1+x)$ for $x\ge 0$. Then  $$g'(x)=\f
1{(1+x)^2}-\f 1{1+x}=-\f x{(1+x)^2}\le 0.$$ Hence $g(x)\le g(0)=0$
and so $F'(x)\le 0$ for $x\ge 0$. Therefore, $F(x)$ and $f(x)$ are
decreasing functions for $x\ge 0$. For $1\le k<n$ we have $\f 1n<\f
1{n-k}$ and so
\begin{align}\label {2.1}
\Big(1+\f 1n\Big)^n\ge \Big(1+\f
1{n-k}\Big)^{n-k}.\tag{2.1}\end{align} Now we prove the inequality
$\b nk\le \f{n^n}{k^k(n-k)^{n-k}}$ by induction on $n$. Since
$k+1\le (k+1)\sls{k+1}k^k$ we see that the inequality holds for
$n=k+1$. Suppose that $\b nk\le \f{n^n}{k^k(n-k)^{n-k}}$. Using
(2.1) we see that
$$\b{n+1}k=\f{n+1}{n+1-k}\b nk\le \f{n+1}{n+1-k}\cdot
\f{n^n}{k^k(n-k)^{n-k}}\le \f{(n+1)^{n+1}}{k^k(n+1-k)^{n+1-k}}.$$
Thus the lemma is proved by induction.

\begin{lem}\label{Lemma 2.2} Let $n\ge 88$. Then
$$\b n{\lceil \frac{n}{3} \rceil }<\f{2^n}{n+1}.$$
\end{lem}

Proof. We prove the result by considering the following three cases:
\begin{itemize}
\item[{\rm Case 1 :}] Suppose that $n=3m$. Taking $n=3m$ and $k=m$ in Lemma 2.1
we get
\begin{align*}
\binom{n}{\lceil \frac{n}{3} \rceil }& =\binom{3m}{m } < \left(
\frac{27}{4}\right)^m = 2^{3m} \cdot  \left( \frac{27}{32}\right)^m
.
\end{align*}
For $m=26$, we have $\left( \frac{27}{32}\right)^m < \frac{1}{3m+1
}$. Since  $\left( \frac{27}{32}\right)^m < \frac{1}{3m+1 }$ implies
that \noindent $\left( \frac{27}{32}\right)^{m+1} < \frac{27}{32
}\cdot \frac{1}{3m+1} < \frac{1}{3(m+1)+1}$, we see that the result
holds for $n\geq 78$.
\item[{\rm Case 2 :}] Suppose that $n=3m+1$. Then  we find that
\begin{align*}
\binom{n}{\lceil \frac{n}{3} \rceil }&=\binom{3m+1}{m+1 } =\frac{3m+1}{m+1} \binom{3m}{m }\\
& < 3\cdot 2^{3m} \cdot  \left( \frac{27}{32}\right)^m
= \frac{3}{2}  \cdot\left( \frac{27}{32}\right)^m \cdot 2^{3m+1}\\
&<  \frac{2^{3m+1}}{3m+2}= \frac{2^n}{n+1},
\end{align*}
for $m\geq 29$ in a similar manner.
\item[{\rm Case 3 :}] Suppose $n=3m+2$. Then we find that
\begin{align*}
\binom{n}{\lceil \frac{n}{3} \rceil }&=\binom{3m+2}{m+1 } =\frac{2}{3}\cdot \binom{3(m+1)}{m+1 }\\
&< \frac{2}{3} \cdot 2^{3m+3} \left( \frac{27}{32}\right)^{m+1} = \frac{4}{3} \cdot 2^{3m+2} \left( \frac{27}{32}\right)^{m+1} \\
&< \frac{2^{3m+2}}{3m+3} = \frac{2^n}{n+1},
\end{align*}
for $m\geq 27$ in a similar manner.
  Thus
 the lemma follows for all integers $n\geq 88$.
\end{itemize}
\noindent{\bf Remark 2.1} Suppose $n\e 0\mod 3$. Taking $m=3$, $p=1$
and replacing $n$ by $n/3$ in [6, Theorem 2.6] we deduce that
$$\b n{\f n3}<\f{3}{2\sqrt{\pi n}}\Big(\f 3{2^{\f 23}}\Big)^n
<0.85\cdot\f{1.89^n}{\sqrt n}.$$ Taking $r=2$ and $m=n/3$ in [4-5]
we derive that
$$\b n{\f n3}>\f {3}{8n}\Big(\f 3{2^{\f 23}}\Big)^n>0.375\cdot
\f{1.88^n}n.$$

\section{\bf Proof of  theorems}
\par {\bf Proof of  Theorem 1.1: } \vskip 3mm
As  $f(n)<\f{n+2}2$, we see that
$$\binom{n+1}{f(n)-1}=\frac{n+1}{n+2-f(n)}\binom{ n}{f(n)-1}
<\frac{n+1}{n+2-f(n)}\cdot
\frac{2^n}{n+1}=\f{2^n}{n+2-f(n)}<\frac{2^{n+1}}{n+2}.$$
 Hence
$f(n+1)\ge f(n)$. On the other hand, as $f(n)\le \f n2$ we have
$$\binom{n+1}{f(n)+1}=\frac{n+1}{f(n)+1}\binom{ n}{f(n)}
>\frac{n+1}{f(n)+1}\cdot\frac{2^n}{n+1}\ge  \frac{2^{n+1}}{n+2}.$$
Thus, $f(n+1)\le f(n)+1$.
\par By (\ref{1.1}),  $L(n)\le [\frac{
n}{2}]\le \frac{ n}{2}$. Thus $L(n)<\frac{ n}{2}+1-\frac{
1}{2n}=n+2-\frac{(n+1)^2}{2n}$ and so
\begin{align*} \binom{n+1}{L(n)-1}=\frac{n+1}{n+2-L(n)}\binom{ n}{L(n)-1}
<\frac{n+1}{n+2-L(n)}\cdot
\frac{2^n}{n}<\frac{2^{n+1}}{n+1}.\end{align*} Hence $L(n+1)\ge
L(n)$. On the other hand, as $L(n)+1<\frac{n}{2}+1+\frac{ 1}{2n}=
\frac{(n+1)^2}{2n}$ we have
$$\binom{n+1}{L(n)+1}=\frac{n+1}{L(n)+1}\binom{ n}{L(n)}
>\frac{n+1}{L(n)+1}\cdot\frac{2^n}n{} > \frac{{n+1}}{(n+1)^2/ (2n)} \cdot\frac{2^n}{n} = \frac{2^{n+1}}{n+1}.$$
Thus, $L(n+1)\le L(n)+1$. This completes the proof.

\vskip 20pt

 {\bf Proof of  Theorem 1.2: } \vskip 3mm
It is sufficient to show that $f(n+1) \leq f(n-1)+1$. This is
equivalent to proving the inequality
$$
\binom{n+1}{f(n-1)+1} > \frac{2^{n+1}}{n+2}.
$$
Now \begin{align*}\binom{n+1}{f(n-1)+1}& =
\frac{n(n+1)}{(f(n-1)+1)(n-f(n-1))} \binom{n-1}{f(n-1)}
\\
 &>\frac{n(n+1)}{(f(n-1)+1)(n-f(n-1))} \cdot \frac{2^{n-1}}{n}\\
& \geq \frac{n(n+1)}{\frac{(n+1)^2}{4}}\cdot
\frac{2^{n-1}}{n}=\f{2^{n+1}}{n+1}> \frac{2^{n+1}
}{n+2}.\end{align*} Thus $f(n+1) \leq f(n-1)+1$. \par By Theorem
1.1, it is sufficient to show that $L(n+1) \leq L(n-1)+1$. This is
equivalent to proving the inequality
$$
\binom{n+1}{L(n-1)+1} > \frac{2^{n+1}}{n+1}.
$$
Now
\begin{align*}\binom{n+1}{L(n-1)+1}& =
\frac{n(n+1)}{(L(n-1)+1)(n-L(n-1))} \binom{n-1}{L(n-1)}
\\
 &>\frac{n(n+1)}{(L(n-1)+1)(n-L(n-1))} \cdot \frac{2^{n-1}}{n-1}\\
& \geq \frac{n(n+1)}{\frac{(n+1)^2}{4}}\cdot \frac{2^{n-1}}{n-1}
> \frac{2^{n+1} }{n+1},
\end{align*}
where we have used the inequality $ab\leq \frac{(a+b)^2}{4}$. The
proof is now complete.
 \vskip 20pt  \noindent {\bf
Proof of Theorem 1.3: } \vskip 3mm
 Since $\b n{L(n)}>\f {2^n}n>\f{2^n}{n+1}$ we see that
$f(n)\le L(n)$. Note that $L(n)\le \f n2<\f{2n^2+5n+1}{3n+1}$. We
then have $\f{L(n)-1}{n+2-L(n)}<\f{2n}{n+1}$. Thus,
$$\b n{L(n)-2}=\f{L(n)-1}{n+2-L(n)}\b n{L(n)-1}\le \f{L(n)-1}{n+2-L(n)}
\cdot\f{2^n}n<\f{2^{n+1}}{n+1}$$
and therefore $f(n)\ge L(n)-1$.

\vskip 20pt \noindent {\bf Proof of  Theorem 1.4: } \vskip 3mm
\noindent As $L(n)\le n/2$ we see that
\begin{align*}\b{n-2}{L(n)-2}&=\f{(L(n)-1)(n+1-L(n))}{n(n-1)}\b
n{L(n)-1}
\\&\le \f{(L(n)-1)(n+1-L(n))}{n(n-1)}\cdot\f{2^n}n
=\f{\f{n^2}4-(\f{n+2}2-L(n))^2}{n(n-1)}\cdot\f{2^n}n
\\&< \f{n^2/4}{n(n-1)}\cdot\f{2^n}n=\f{2^{n-2}}{n-1}.\end{align*}
Thus, $f(n-2)\ge L(n)-1$. By Theorem 1.1 and Corollary 1.1,
$f(n-2)\le f(n-1)\le f(n)=L(n)-1$. Thus $f(n-2)=f(n-1)=
f(n)=L(n)-1$.
\par As
\begin{align*} \b{n+2}{L(n)}&=\f{(n+1)(n+2)}{L(n)(n+2-L(n))}\b
n{L(n)-1}
>\f{(n+1)(n+2)}{L(n)(n+2-L(n))}\cdot \f{2^n}{n+1}
\\&\ge \f{(n+1)(n+2)}{(n+2)^2/4}\cdot\f{2^n}{n+1}=\f{2^{n+2}}{n+2},
\end{align*}
we see that $L(n+2)\le L(n)$. Since $L(n+2)\ge L(n+1)\ge L(n)$ by
Theorem 1.1, we get $L(n)=L(n+1)=L(n+2)=f(n)+1$. This proves the
theorem.

\vskip 20pt \noindent {\bf Proof of  Theorem 1.5: } \vskip 3mm
\noindent By the definition of $f(n)$, we only need to prove that
$$\b{n+3}{f(n)}<\f{2^{n+3}}{n+4}.$$ It is clear that
\begin{align*}
\b{n+3}{f(n)}&=\f{(n+3)(n+2)(n+1)}{f(n)(n+3-f(n))(n+2-f(n))}\b
n{f(n)-1}\\&\le \f{(n+1)(n+2)(n+3)}{f(n)(n+3-f(n))(n+2-f(n))} \cdot
\f {2^n}{n+1}.\end{align*} For $1\le x\le\f n2$ let
$$F(x)=x(n+3-x)(n+2-x)=x^3-(2n+5)x^2+(n+2)(n+3)x.$$ Then
$$F'(x)=3x^2-2(2n+5)x+(n+2)(n+3)=3(x-x_1)(x-x_2),$$
where
$$x_1=\f{2n+5-\sqrt{(n+\f 52)^2+\f
34}}3<\f{2n+5-(n+\f 52)}3<\lceil\f n3\rceil+1$$ and
$$x_2=\f{2n+5+\sqrt{(n+\f 52)^2+\f 34}}3\ge n.$$ For $\lceil \f n3\rceil+1\le x\le\f n2$
we have $x_1<x\le \f n2<x_2$ and so $F'(x)<0$. Thus, $F(x)$ is an
decreasing function and hence
$$F(x)\ge F(\f n2)=\f n2(n+3-\f n2)(n+2-\f n2)=\f 18n(n+4)(n+6).$$
\par By Lemma 2.2,  $f(n)\ge \lceil\f n3\rceil+1$ for $n\ge 88$. We then get
$$f(n)(n+3-f(n))(n+2-f(n))=F(f(n))\ge \f 18n(n+4)(n+6).$$
Now, from the above we deduce that
$$ \b{n+3}{f(n)}\le \f{(n+1)(n+2)(n+3)}{f(n)(n+3-f(n))(n+2-f(n))} \cdot
\f {2^n}{n+1}\le
\f{(n+1)(n+2)(n+3)}{n(n+4)(n+6)}\cdot\f{2^{n+3}}{n+1}.$$ As $n\ge
88$, we have $(n-\f 52)^2>47> 40+\f {25}4$ and so $n(n-\f
52)^2-(39+\f{25}4)n-18>n-18>0$. That is, $n^3-5n^2-39n-18>0$. Thus,
$$n(4n^2+13n-6)-3(n+1)(n+2)(n+3)=n^3-5n^2-39n-18>0$$ and so
$$\f{n(n+4)(n+6)-(n+1)(n+2)(n+3)}{(n+1)(n+2)(n+3)}
=\f{4n^2+13n-6}{(n+1)(n+2)(n+3)}
>\f{n+3-n}n.$$
Therefore,
$$\f{n(n+4)(n+6)}{(n+1)(n+2)(n+3)}>\f{n+3}n\qtq{and so}
\f{(n+1)(n+2)(n+3)}{n(n+4)(n+6)}<\f n{n+3}.$$ Hence, from the above
we deduce that
$$ \b{n+3}{f(n)}\le
\f{(n+1)(n+2)(n+3)}{n(n+4)(n+6)}\cdot\f{2^{n+3}}{n+1}<\f n{n+1}
\cdot\f{2^{n+3}}{n+3} <\f{2^{n+3}}{n+4}.$$ This proves the theorem.
\vskip 20pt

\noindent\textbf{Acknowledgements}: The second and third authors
wish to thank the National Institute for Mathematical Sciences
(NIMS), Daejeon, Republic of Korea for its warm hospitality and
generous support. The corresponding author Zhi-Hong Sun is supported
by the Natural Sciences Foundation of China (grant no. 11371163).

\end{document}